# Decision Making: Lexicographical Procedure.

V. Zhukovin, N. Chkhikvadze, Z. Alimbarashvili,

**Abstract**
It is introduced the using of generation lexicographical procedure for multicriteria.  decision-making problems.

## 1. Introduction

The lexicographical procedure widely spread in the decision making problems. It is used in the various sphere of the human action. for example:

1. The words are lexicographically regulated in the dictionary.
2. The base-ten system of the number recording is defined by lexicographical principle.
3. The militaries use the lexicographical principle in the problem of the objective allocation.
4. The American scientists experimentally proved, that the humans use the lexicographical principle for decision making.
5. It is experimentally proved, that the dolphins use the lexicographical principle in the problem of the choice. ( T. Zorikov, institute of cybernetics Georgian academy of sciences.)

Lexicographical procedure of choice holds special space in the theory of decision making. The computer realization of lexicographical procedure will be successfully use in various practical problems.

Analysis of lexicographical procedure is conducted for a long time. One of interesting is work of .Podinovcki (Moscow) B. Zhukovin (Institute of cybernetics Georgian academy of sciences ) use lexicographical procedure in unclear problems of decision making. But as though full it was not investigated this scientific direction, always there are unexplored problems, both in the theory  and in computer realization. It causes our present interest to this procedure. Its basic advantage this reduction of time of decision-making, and also reduction of volume of the information necessary for decision-making. It is usually effectively used in systems, where small operative memory.

## 2. Primary formulation.

From the beginning we formulate general lexicographical principle and then use it in various procedure of choice. For the obviousness we imagine that we have some conditional system with main memory $\upsilon$ (volume). We solve some problem of the decision making and the information with volume $V$ ( all the concrete information which is connected with problem) may be used for it. Note, that the introduced values are not numerical characteristic, though like by implication. Assume $\upsilon < V$ at the same time we can not use all the information for one step (cycle). We can use only part of  the information , which may be located in $\upsilon$. So all the volume of the information $V$ are separated into blocks: $V_1, V_2, ..., V_t$ . These separating must satisfy on the following condition:

1. $V_j \leq \upsilon$ for all $j = 1 \dot{-} t$

2. $$\bigcup_{j=1}^{t} V_j = V$$
3. The blocks must be whole ( for example the part of the table can not place into one block and other into different block).
4. The blocks must be ordering by importance of the information.

Really there are not all the conditions, but the law of lexicographical idea consists in these 4 conditions. The matter is that often at decision making in concrete areas the choice can be made using only a part of the information and neglect all other. Thus quality the decision does not suffer.

Now may be formulated itself lexicographical principle:
1. Let all the condition fulfilled – the problem get ready to lexicography.
2. Decision making performs in sequence by steps. (It may be not more $t$)
3. The block of the information $V_j$ is used at $j$ - step. Direction of the work is: $1 \to t$
4. The result: The choice made, or is not indefinite - full information $V$ is not sufficient for the choice. There is needed side information.

**?** Clearly that the time of decision making and the information using during this time reduce because the choice is realized earlier than at $t$ step.

### 3. Multicriteria lexicographical procedure.

Multicriteria decision making problem has been mathematical formalized with L. Zade in 1960. After that lexicographical procedure of choice has mathematical base for its description, development and generation new procedures.

Multicriteria problem of decision making is presented in vector interpretation in this way:
$$D = \langle X, K \rangle, \qquad (1)$$
where, $X$ - is finite set of competitive alternatives $x_i \subset X$, $i = 1. \div n$. Let introduce the set of ordered pairs $E$ with the elements $(x_i, x_k) \in E$, $x_i, x_k \in X$.
$$K = \{K_1(x_i), K_2(x_i), ..., K_j(x_i), ..., K_m(x_i)\} \qquad (2)$$

This is vector criterion of efficiency; where $K_j(x_i)$, $j = 1 \div m$, are scalar functions defined on the set $X$. We can suppose, that all $K_j(x_i)$ are type of " win". Besides they are defined with the content and although each of them has its scale of measuring the type of scale is identity.

Since we always can trade places the components at vector criterion $K$ of importance, we could say that $K_j(x_i)$ criterions are ordered by the relevance. The order concurs with the number's order, in other words most important criterion is $j = 1$, next $j = 2$ and so distances $j = m$. It is implied that this order is lexicographical. We explain this affirmation below. Now we formulate lexicographical structure on the set $E$. For it is necessary to define binary lexicographical preference relations $R_{lex}$. It defines in the following way. Let consider the decisions pairs $(x_i, x_l) \in E$.

$R_{lex}$ : It is sad, that the decisions $x_i \in X$, are lexicographical preference than $x_l \in X$, if one of following conditions fulfils:
1) $K_1(x_i) > K_1(x_l)$ , or

2) $K_1(x_i) = K_1(x_l)$ and $K_2(x_i) > K_2(x_l)$, or  (3)

-----------------

j) $K_1(x_i) = K_1(x_l)$ and ...and $K_{j-1}(x_i) = K_{j-1}(x_l)$ and $K_j(x_i) > K_j(x_l)$ or

-----------------

n) $K_1(x_i) = K_1(x_l)$ and ... and $K_{n-1}(x_i) = K_{n-1}(x_l)$ and $K_n(x_i) > K_n(x_l)$

The stop happens on the row where the inequality is fulfilled. If the inequality is not executed on $n$ row too, i.e. $K_n(x_i) = K_n(x_l)$, then this pair is lexicographically equivalent decision. Briefly this procedure is designating in that way:

$$x_i \overset{lex}{\underset{\sim}{\succ}} x_l \qquad (4)$$

Let explore characteristics of $R_{lex}$ In the first place it is **linked,** i.e. all decision pairs from $E$ are congruous with this relation. It is also **asymmetrical** and **transitive**.

Superiority relation $R_{lex}$ is equivalence too. Such superiority relation is linear order i.e. $R_{lex}$ is linear order. Now we shall apply all arsenal of result, theorems, which are available at present in the theories of decision making established us and other scientists, to research of the presented lexicographical procedure (3). And here the first. It is proved that the linear order has nonempty set Pareto (its kernel) and this set contains one decision, if it some they are equivalent. So procedure (3) always give out one decision. Whether it is the best (utopian) by M. Saluqvadze ( A. Eliashvili institute of control systems Georgian academy of sciences), but other mathematical representations of procedure (3) exist. We investigate one of them (a kind of linear convolution of $K$ vector criterion) below.

### 4. Lexicographical coefficients of importance.

Criterion $K_j(x_i)$ of efficiency are ordered on importance according to numbers: the less number the more important criteria for a choice. This fact can be fixed, having attributed criteria coefficients of importance: $\lambda_j \to K_j(x_i)$, $j = 1 \div m$. In this case in the theory of decision making linear convolution of vector criteria of efficiency $K$ is considered:

$$L(x_i) = \sum_{j=1}^{m} \lambda_j * K_j(x_i) \qquad (4)$$

This convolution possesses many very good properties which we shall consider later, in the same section. The unique requirement to importance coefficients are $\lambda_j \geq 0$ for all $j = 1 \div m$. Sometimes use also a condition $\sum_{j=1}^{m} \lambda_j = 1$, but it not definitely. With this the condition is executed:

if $L(x_i) \geq L(x_\ell)$ then $x_i \underset{\sim}{\succ} x_\ell$, and  (5)

on the contrary if $x_i \underset{\sim}{\succ} x_\ell$ then $L(x_i) \geq L(x_\ell)$.

Clear that if $\lambda_j$ are numbers and $K_j(x_i)$ scalar function, then $L(x_i)$ scalar function definite on $X$. The freedom of a choice $\lambda_j, j = 1 \div n$ in this case us is not arranged.

Self-wiled collection of important coefficients $\Lambda = \{\lambda_j\}_1^m$ disturbs lexicographical conditions representing by the formula (3).Whether there is a question is there such collection of importance coefficients which would not break, and kept a condition lexicographic (3)? Yes, it is and not only one but whole class. Therefore we designate lexicographical coefficients of importance with $\lambda_j(\ell ex)$, $j = 1 \div m$ and all collection of coefficients with

$$\{\lambda_j(\ell ex)\}_1^m = \Lambda(\ell ex). \tag{6}$$

When linear convolution writes so for lexicographical choice (basic formula):

$$L_{\ell ex}(x_i) = \sum_{j=1}^m \lambda_j(\ell ex) * K_j(x_i). \tag{7}$$

There are worked special procedures for the accounting of importance coefficients $\lambda_j(\ell ex)$, $j = 1 \div m$, by several scientists and by us too. These procedures do not coincide and use different scales of measuring. Now we describe basic property of linear convolution:

a) $L_{\ell ex}$ is a function of usefulness definite on the $X$. It assigns linear order according to rules:

$$L_{\ell ex}(x_i) \geq L_{\ell ex}(x_\ell) \leftrightarrow x_i \underset{\sim}{\overset{lex}{\succ}} x_\ell \tag{8}$$

b) This order is adjusted with the order presented by formula (3), this is:
   if $x_i \underset{\sim}{\succ} x_\ell$ correspondence with formula (3), that

$$L(x_i) \geq L(x_l) \tag{9}$$

and on the contrary. The formulas (8) and (9) are identical.

c) Any new presentation must be adjusted with (3) and it means that the presentation will be adjusted with (7) too.

d) Pareto set $X_\Pi(\ell ex)$ contains one decision (kernel of linear orders) . This best decision is single and does not depend for the view of the presentation.

e) This best element will be find, if we decide simple optimization problem:

$$x^* = \max_{x_i \in X} L_{\ell ex}(x_i) \tag{10}$$

where $x^* \in X$ and it is best decision in lexicographical procedure of choice (3).

f) Now about scale of measure. All known us collection of lexicographic coefficients are presented by ordering scale:

$$Ш = \langle K_j, \Lambda_{\ell ex}, \Phi(K_j) \rangle, \tag{11}$$

with this $\Phi$ is the class of permitting transformation. It is monotone function (increasing or decreasing). Usually the decision $x^*$ is instable in the ordering scale i.e. the expression (10) does not invariant concerning monotone transformation of the coefficients. But this common rule breaks in lexicographical procedure The decision $x^*$ is stability concerning the transformation. The expression (10) is invariant in ordering scale. It is related to this max is searching with one criterion and which is defined by compare pair of decision.

g) Cleverly, lexicographical choice link to superiority degree: $Z(x_i, x_\ell) = j_0$
   ( it is number of line where is executed the choice)

Let analyze detail the criterions $K_j$. It is scalar function definite on $X$. Let consider these coefficients. Every criterion has own name (energy, length, expenditure), which represents any property of the estimating object $x_i \in X$. The name enables us the group of criterion $K_j, j = 1 \div m$, will be ordered with importance.

Every number criterion has the scale of estimating. There is example pic.1.

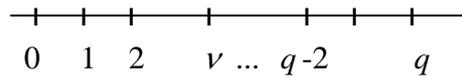

0   1   2   $\nu$ ...  q-2   q

pic.1.

Usually the scale is beginning from the zero, but it is possible beginning from either of meaning (for example $p$). Let consider the criterion of "win". It means that $\nu > \nu - 1$ where $\nu$ is indicator of scale rank. This dependence is transitive and linked. Now define scale diapason of criterion $K_j$.

$$d_j = \max \nu(j) - \min \nu(j), \qquad (12)$$

if the scale is beginning from zero $d_j = \max \nu(j) = q$, where $q$ is quantity of scale rank.

Let define the conditions of lexicographical group criterions $K_j, j = 1 \div m$:

**Affirmation 1.**

Assume the importance of criterion is increasing from the right to the left and its number of position from the left to the right (formulae 3), then if the condition $\min \nu(j) \succ \max \nu(j+1)$ is fulfilled for all $j = 1 \div m - 1$, then criterions group $K_j, j = 1 \div m$ is lexicographic ordered.

Notice:
1. Zero takes not in part of the minimum defining (only in the affirmation)
2. The mark $\succ$ signifies the "importance"
3. This dependence is linked and transitivity.

Now we can assign the formula for the accounting lexicographical coefficients of importance:

$$\lambda_j(\ell ex) = d_j^{j-1}, \; j = 1 \div m \qquad (13)$$

May be the scales of the criterion are diverse (for example: the speed m/h, energy k/h). In this case the formulae (3) works, but $L(\ell ex)$ loses the sense. One of the method the forming uniformity system of the criterion is the rationing. Instead of $K_j$, include following criterions:

$$\gamma_j(x_i) = \frac{K_j(x_i)}{\max K_j} * a \quad \text{where } a > 0 \qquad (14)$$

This gives the freedom for the variation. ( for example $a = 10$ scale of Saaty).

### 6. Non- linked lexicographical procedure of decision making

The new type of binary relation allows easily pass to indistinct lexicographical relations of preference. let replace number criterions of effectiveness $K_j$ by binary relations of preference $R_j$ in the formulae (3). These are ordered with the importance. It makes the expert. All previous definitions and results lose effect for this class. Full research must be make for beginning (for definition of lexicographical non-linked binary relations of preference). Either $R$ corresponds $R^{-1}, R^e, R^s, R^N$ are incomparable pairs and because of it we have non-linked relation.

**Definition:**

The pairs $(x_i, x_\ell) \in E$ are comparing, if

$$(x_i, x_\ell) \in R_1^s, \text{ or}$$
$$(x_i, x_\ell) \in R_1^s, \text{ and } (x_i, x_\ell) \in R_2^s, \text{ or} \qquad (15)$$
$$(x_i, x_\ell) \in R_1^s, \ (x_i, x_\ell) \in R_2^s, \ \ldots \text{ and } \ (x_i, x_\ell) \in R_m^s, \text{ then}$$

they say $x_i$ is strongly preference $x_\ell \in X$. If $(x_i, x_\ell) \in R_m^\ell$, then $x_i$ and $x_\ell$ are equivalence. If $(x_i, x_\ell) \in R_j^N$, is appeared on any step of comparing process then we will tell that decision pair is not comparing by lexicographical. In this case lexicographical relation of importance (which is presented by formulae (15)) is non-linked.

Let form two axioms.

$A_1$: If one of two $x_i$ and $x_q$ equivalence decisions, $x_i$ is importance than any decision $x_\ell \in X$, i.e. $(x_i, x_\ell) \in R_j^s$, $j = 1 \div m$, then second decision $x_q$ is importance than $x_\ell$, i.e. $(x_q, x_\ell) \in R_j^s$,

$A_2$: (opposite) If any decision $x_\ell \in X$ is importance than one of two $x_i$ and $x_q$ equivalence decisions for example $x_i$ i.e. $(x_\ell, x_i) \in R_j^s$, $j = 1 \div m$, then it is importance than second decision $x_q$ i.e. $(x_\ell, x_q) \in R_j^s$.

**Affirmation 2.**

If $R_j$, $j = 1 \div m$, are non-linked, transitive and fulfilled the axioms A1 and A2, then lexicographical relation of importance definition by formulae (15) is transitive and non-linked too.

In this case we tell about of quasi- ordering. In general, quasi-ordering may be presented by the function preserving order (F P O), but how is possibility formed for this lexicographical relation of importance is not clear. For now we consider two variants: through characteristic function of relation and with the using superiority degree.(Zhukovin, Institute of Cybernetics, Tbilisi). This problem will be required subsequent researching, also the finding conditions for the examine the groups of relation for the lexicographic.

### decision-making: Indistinct lexicographic procedure

First of all this procedure should be defined. Some definitions will be necessary for this purpose from the theory indistinct sets (L.Zade) for us.

$X$ - on former set of competitive decisions. Let it certainly - for practical problems this requirement is natural.

$T = X \times X$ set of all ordered pairs decisions. If $X$ contains $n$ competitive decisions: $x_i \in X$ ; $i = 1 \div n$ pairs will be $n^2$, including pairs $(x_i, x_i) \in E$. On set $E$ we shall define the indistinct binary attitude:

$$M = \langle E; \mu(x_i; x_\ell) \rangle \qquad (16)$$

where $\mu(x_i; x_\ell)$ there is the function of the accessory accepting values from the closed interval [0;1].

As $E$ during all procedure does not vary, a basic element $M$ is function of an accessory. We already know, that to any binary attitude including careless correspond:

a) $M^{-1} = \langle E; \mu^{-1}(x_i; x_\ell) \rangle$ - the return attitude,

b) $M^\ell = \langle E; \mu^\ell(x_i; x_\ell) \rangle$ - attitude equivalences, (17)

c) $M^s = \langle E; \mu^s(x_i; x_\ell) \rangle$ - attitude superiority.

Thus:

a) $\mu^{-1}(x_i; x_\ell) = \mu(x_\ell; x_i)$ ;

b) $\mu^\ell(x_i; x_\ell) = \min\{\mu(x_i; x_\ell); \mu(x_\ell; x_i)\}$; (18)

c) $\Delta(x_i x_{\ell_i}) = \mu(x_i; x_\ell) - \mu(x_\ell; x_i)$

$$\mu^s(x_i; x_\ell) = \begin{cases} \Delta(x_i; x_\ell) , & \text{if } \Delta(x_i; x_\ell) \geq 0 \\ 0 , & \text{if } \Delta(x_i; x_\ell) \leq 0 \end{cases}$$

Now we can define indistinct lexicographic procedure: (FLEX) unlike precise lexicographic procedure (LEX).

We compare two competitive decisions: $x_i \in X$ and $x_\ell \in X$, or (that –also) is worked with pair $(x_i, x_\ell) \in E$. Comparison goes by many criteria. Let them will be $m$. In this case as these criteria indistinct attitudes of preference (AP), or (that - too most) the indistinct vector attitude of preference (VAP) are used. These AP are ordered on importance lexicography: the less its number, the it is more important. Comparison is carried out on steps. Number of a step coincides with number AP which works on the given step. Now we shall result procedure of comparison:

1) $\mu_1^s(x_i; x_\ell) > 0$ or

2) $\mu_1^s(x_i; x_\ell) = \mu_1^s(x_\ell; x_i) = 0$ and $\mu_2^s(x_i; x_\ell) > 0$ or

- - - - - - - - - - - (19)

j) $\mu_q^s(x_i; x_\ell) = \mu_q^s(x_\ell; x_i) = 0$, $q = 1 \div (j-1)$ and $\mu_j^s(x_i; x_\ell) > 0$ , $j = 1 \div m$

At j=1 and at j=2 it is received two top lines. And at j =m it is received last line. Thus into this definition enters $m$ lines. It describes indistinct strict lexicographic preference of the decision $x_i$ above the decision $x_\ell$. Thus if comparison was carried out at line $j_0$.

$$\mu_{\ell tx}^s(x_i; x_\ell) = \mu_{j_0}^s(x_\ell; x_i) \qquad (20)$$

If comparison was not carried out even on last step $x_i$ and $x_\ell$ lexicography equal(are equivalent). In this case write $\mu_{\ell ex}^\ell(x_i; x_\ell)$.

We shall consider some properties of the indistinct lexicographic attitude of preference: A) It is transitive and coherently. It means, that it is linear order. Its set Pareto not empty. B) It can be presented in the form of generalize linear criterion.

$$\mu_{\ell tx}^{s}(x_i; x_\ell) = \sum_{j=1}^{m} \lambda_{\ell ex}^{j}\, \mu_{j}^{s}(x_i; x_\ell) \tag{21}$$

Lexicographic coefficients $\lambda_{\ell ex}^{j}$ it is not obviously possible to calculate by means of above resulted [1] algorithms as the basis at all $\mu^s(x_i; x_\ell)$ is equal 1. We shall resort to an artificial method. An interval [0; 1] on which defined these functions of an accessory, we shall break into 10 intervals in length 0,1. We shall make That the basis of the factors, equal 10. And then by analogy to a decimal notation we shall receive for factors the formula:

$$\lambda_{\ell ex}^{j} = 10^{m-j} \tag{22}$$

The greatest will be $\lambda_{\ell ex}^{1} = 10^{m-1}$, it corresponds the most to the important attitude of preference: $\mu_1^s(x_i; x_\ell)$. Most small will be: $\lambda_{\ell ex}^{m} = 10^0 = 1$, it corresponds to the least important attitude of preference.

**About scales.** The lexicographic vector contains $m$ a component, each of which is the linear order. It means, that they are set in a serial scale. In a precise case for $K_j(x_i)$ and $R_j$ it is practically obvious. And in an indistinct case for $\mu_j(x_i; x_\ell)$ it is possible to enter the following size:

$$\mu_j(x_i) = \frac{1}{n} \sum_{x_\ell \in X} \mu_j(x_i; x_\ell) \tag{23}$$

It is one-dimensional function of an accessory(a belonging). It can be considered(examined) as function of utility. Let there is lexicographic vector $L_1$ with lexicographical ordered components $\gamma_j$, $j = 1 \div m$ (the less number, the more important a component). It forms on $X$ lexicographic structure $StrLex_1$. Everyone a component is set in a serial scale. We apply the resolved(allowed) transformation of a scale, same for all $m$ a component. Let's receive new lexicographic vector $L_2$ with components $r_j$, $j = 1 \div m$ to which there will correspond(meet) lexicographic structure of $StrLex_2$. The resolved(allowed) transformation for serial scales is strictly growing, or strictly decreasing function.

**The theorem.**    $StrLex_1 = StrLex_2$

Equality means, that structures coincide.

**The conclusion.** Theoretical research of lexicographic procedures of comparison and a choice has allowed us to develop the first variant of dialogue system of formation of any dictionaries.[3] Thus on a computer the most labour-consuming and tiresome area of work is realized. This ordering of words. Most likely, such programs exist, but we is constructed on other basis and differs variety of decided(solved) problems(tasks). Besides we are going to expand her(it) and on other areas of human activity where it is applied lexicography.